\RequirePackage{ifthen}
\newboolean{esa}
\setboolean{esa}{false}

\ifthenelse {\boolean{esa}}
{
	
	\documentclass[a4paper,UKenglish,cleveref, autoref,  thm-restate]{lipics-v2021}
	
} {
	\documentclass[a4paper,12pt]{article}
	\usepackage[a4paper, left=3cm, right=3cm, top=2.5cm, bottom=2.5cm]{geometry}

}
\usepackage{amsmath}
\usepackage{graphicx}
\usepackage{psfrag}
\usepackage{amssymb}
\usepackage{amsmath}
\usepackage{verbatim,booktabs}
\usepackage{epsfig,url}
\usepackage{float}
\usepackage{mathrsfs}
\usepackage{caption}
\usepackage{tabularx}

\usepackage[usenames,dvipsnames]{pstricks}
\usepackage{epsfig}
\usepackage{bbm}
\usepackage{comment}
\usepackage{multirow}
\usepackage{array}
\usepackage{algpseudocode} % http://ctan.org/pkg/algorithmicx
\usepackage{algorithm}% http://ctan.org/pkg/algorithm
\usepackage{xspace}% http://ctan.org/pkg/xspace
\usepackage{hyperref}
\usepackage{cleveref}
\usepackage{todonotes}%

\makeatletter
\let\cl@chapter\undefined
\makeatletter

\algnewcommand{\algorithmicgoto}{\textbf{go to}}%
\algnewcommand{\Goto}[1]{\algorithmicgoto~\ref{#1}}%

\def\ie{{i.e.,} }

\newcommand{\fnm}[1]{#1}%
\newcommand{\sur}[1]{#1}%

\newcolumntype{L}{>{$}l<{$}}
\newcolumntype{C}{>{$}c<{$}}
\newcolumntype{R}{>{$}r<{$}}

\newcommand{\bR}{\mathbb{R}}
\newcommand{\bZ}{\mathbb{Z}}

\newcommand{\norm}[1]{{\lVert#1\rVert}}

\newcommand{\relx}[1]{\tilde{#1}}

\DeclareMathOperator{\Diag}{Diag}

\DeclareMathOperator{\bd}{bdd}

\newcommand{\floor}[1]{\lfloor #1 \rfloor}

\Crefname{chapter}{Chap.}{Chaps.}
\Crefname{section}{Sect.}{Sects.}
\Crefname{proposition}{Prop.}{Props.}
\Crefname{theorem}{Thm.}{Thms.}
\Crefname{definition}{Defn.}{Defns.}
\Crefname{corollary}{Cor.}{Cors.}
\Crefname{figure}{Fig.}{Figs.}
\Crefname{observation}{Obs}{Obss.}

\newcommand{\win}[1]{\textbf{\textcolor{blue}{#1}}}
\newcommand{\loss}[1]{\textit{\textcolor{red}{#1}}}

\title{Sparsity-driven Aggregation of Mixed Integer Programs} %TODO Please add

\ifthenelse {\boolean{esa}}
{
	
	\author{Liding Xu\footnote{corresponding author}}{Zuse Institute Berlin, Germany}{lidingxu.ac@gmail.com}{https://orcid.org/0000-0002-0286-1109}{}
	
	\author{Gioni Mexi}{Zuse Institute Berlin, Germany}{mexi@zib.de}{https://orcid.org/0000-0003-0964-9802}{}
	
	\author{Ksenia Bestuzheva}{Zuse Institute Berlin, Germany}{bestuzheva@zib.de}{https://orcid.org/0000-0002-7018-7099}{}
	
	\authorrunning{Liding Xu et al.} %TODO mandatory. First: Use abbreviated first/middle names. Second (only in severe cases): Use first author plus 'et al.'
	
	\Copyright{CC-BY} %TODO mandatory, please use full first names. LIPIcs license is "CC-BY";  http://creativecommons.org/licenses/by/3.0/

	
	 %TODO mandatory; please add comma-separated list of keywords
	
	\category{} %optional, e.g. invited paper
	
	\relatedversion{} %optional, e.g. full version hosted on arXiv, HAL, or other respository/website
	\funding{The work for this article has been conducted in the Research Campus MODAL funded by the German Federal Ministry of Education and Research (BMBF) (fund numbers 05M14ZAM, 05M20ZBM).}%optional, to capture a funding statement, which applies to all authors. Please enter author specific funding statements as fifth argument of the \author macro.
	
	\acknowledgements{We want to thank Felipe Serrano for his insightful comments on this paper.}%optional
	
	\EventEditors{John Q. Open and Joan R. Access}
	\EventNoEds{2}
	\EventLongTitle{42nd Conference on Very Important Topics (CVIT 2016)}
	\EventShortTitle{CVIT 2016}
	\EventAcronym{CVIT}
	\EventYear{2016}
	\EventDate{December 24--27, 2016}
	\EventLocation{Little Whinging, United Kingdom}
	\EventLogo{}
	\SeriesVolume{42}
	\ArticleNo{23}
}
{
	\newtheorem{example}{Example}

	\newcommand{\orgname}[1]{#1}%

	\title{Sparsity-driven Aggregation of Mixed Integer Programs}

	\date{\today}
	
	\author{\fnm{Liding} \sur{Xu} \thanks{  \orgname{Zuse Institute Berlin, Germany}.
			E-mail: {\tt lidingxu.ac@gmail.com, mexi@zib.de,bestuzheva@zib.de}}
		\and
		\fnm{Gioni} \sur{Mexi} \footnotemark[1]
		\and
		\fnm{Ksenia} \sur{Bestuzheva} \footnotemark[1]
	}
	
}
\begin{document}
	
	\maketitle

	\begin{abstract}
		Cutting planes are crucial for the performance of branch-and-cut algorithms for solving mixed-integer programming (MIP) problems, and linear row aggregation has been successfully applied to better leverage the potential of several major families of MIP cutting planes. This paper formulates the problem of finding good quality aggregations as an $\ell_0$-norm minimization problem and employs a combination of the lasso method and iterative reweighting to efficiently find sparse solutions corresponding to good aggregations. A comparative analysis of the proposed algorithm and the state-of-the-art greedy heuristic approach is presented, showing that the greedy heuristic implements a stepwise selection algorithm for the $\ell_0$-norm minimization problem. Further, we present an example where our approach succeeds, whereas the standard heuristic fails to find an aggregation with desired properties. The algorithm is implemented within the constraint integer programming solver SCIP, and computational experiments on the MIPLIB 2017 benchmark show that although the algorithm leads to slowdowns on relatively ``easier'' instances, our aggregation approach decreases the mean running time on a subset of challenging instances and leads to smaller branch-and-bound trees.
	\end{abstract}

	\section{Introduction}
	Mixed integer programs (MIPs) are an important class of mathematical programs combining continuous and discrete decision variables with linear objective functions and constraints.
	The branch-and-cut algorithm is the backbone technique used in many MIP solvers, e.g., SCIP \cite{bolusani2024scipoptimizationsuite90}.
	A critical component of the branch-and-cut algorithm's success is the use of cutting planes.
	There are rich classes of valid inequalities for the mixed-integer linear program (MILP):
	\begin{equation}
		\label{eq.milp}
		\min \{c^\top x \,|\, x  \in \bR^p \times \bZ^{q}, Ax - b \le 0\},
	\end{equation}
	where $\bR$ denotes the set of reals,  $\bZ$ denotes the set  of integers, $A \in \bR^{m \times (p+q)}$, and $b \in \bR^m$.

	In this paper, we investigate the common aggregation procedure \cite{marchand2001aggregation} used in the generation of valid inequalities for MILPs.
	We formulate the aggregation problem as a cardinality optimization problem, and we exploit the sparse optimization technique called least absolute shrinkage and selection operator (lasso) \cite{tibshirani1996regression}.
	Our primary contribution is the development of a sparsity-driven and optimization-based aggregation algorithm tailored for MILPs. This allows us to derive more effective complemented Mixed-Integer-Rounding (c-MIR) cuts \cite{marchand19990,nemhauser1990recursive} than the previous heuristic \cite{marchand2001aggregation}.

	We next review the aggregation procedure for MILPs.
	Aggregation constructs a single-row surrogate relaxation, defined by the so-called base inequality, in the form: $ \min \{c^\top x \,|\, x \in \bR^p \times \bZ^{q}, \lambda^\top(Ax - b) \le 0\}$, where $\lambda$ is a vector of non-negative factors. This surrogate relaxation simplifies the structure of the original problem, facilitating the derivation of valid inequalities. An illustrative example is the Chvátal-Gomory (CG) cut \cite{chvatal1973edmonds,fischetti2007optimizing}. For pure integer linear programs (ILPs), CG cuts take the form of $\floor{\lambda^\top A} x - \floor{\lambda^\top b} \le 0$ with $\lambda \in \bR^m_+$, and these cuts are a specific type of split inequality  \cite{cook1990chvatal}, inferred from the aggregated single-row surrogate relaxation. Subclasses of CG cuts include zero-half cuts \cite{caprara19960}, where factors of $\lambda$ are restricted to  zero or one-half, and mod-$k$ cuts \cite{caprara2000separation}, with factors of $\lambda$ in $\{0, 1/k, \dots, k-1/k\}$. For general MILPs, projected CG cuts can be used.  Similar to the principles of  Benders decomposition \cite{bnnobrs1962partitioning}, Bonami et al. \cite{bonami2008projected} proposed projecting out continuous variables from MILPs to reduce them to ILPs, \ie the first $p$ elements of $\lambda^\top A$ corresponding integer variables are zeros, then CG cuts can be applied effectively.

	In addition to split cuts, Mixed-Integer-Rounding (MIR) cuts \cite{nemhauser1990recursive} form another fundamental class of valid inequalities, particularly relevant for single-row MIP models involving continuous variables. While MIR cuts are mathematically equivalent to Gomory’s mixed-integer cuts \cite{gomory1960algorithm}, split cuts, and a special case of disjunctive cuts \cite{balas1979disjunctive}, their separation procedures differ significantly, particularly when combined with aggregation, resulting in distinct implementations in optimization solvers.
	Aggregation is also used for generating complemented MIR (c-MIR) cuts, as introduced in the heuristic by Marchand and Wolsey (MW) \cite{marchand19990,marchand2001aggregation}. The MW heuristic has been shown to be effective in the computational studies realized in SCIP \cite{achterberg2007constraint,wolter2006implementation}, MOPS \cite{christophel2009separation}, and CPLEX \cite{gonccalves2005implementation}.
	The MW heuristic uses a stepwise selection process to create a sparse base inequality that has as few continuous variables as possible, \ie some continuous variables are greedily projected out. %This results in c-MIR cuts that are sparser than CG cuts.\todo{Are the cuts themselves actually sparser, though? All procedures I've seen (including the detailed descriptions of algorithms by Wolter) introduce one variable for the sum of all nonnegative continuous slacks, so, unless actual code does something differently, the cuts always have only one continuous variable in the end. But I also see that you are here specifically comparing to CG cuts. May be some citation or an explanation to back this statement?}
	For MILP instances with network structures, such as network design problems, a tailored heuristic \cite{achterberg2010mcf} has been proposed to find aggregations corresponding to sparse cut sets in networks, and the resulting c-MIR cuts are the well-known flow-cutset inequalities \cite{atamturk2002capacitated,bienstock1996capacitated}.
	Note that the technique of projecting out continuous variables via row aggregations before applying c-MIR cuts is also crucial for enabling other components, such as cut-based conflict analysis \cite{mexi2024cut}, to generate conflict constraints in the presence of continuous variables.

	Our innovation is centered on the use of sparse optimization to find aggregations.  The corresponding sparse optimization problem can be formulated as a cardinality optimization problem \cite{tillmann2024cardinality} involving the nonconvex and non-differentiable $\ell_0$-norm, which is a hard combinatorial optimization problem. As a practical alternative to $\ell_0$-norm minimization, we employ approximations using sparsity-inducing convex norms. In the case of the $\ell_1$-norm,  the approximation problem is lasso \cite{tibshirani1996regression}, which can be solved efficiently by linear programming (LP). Under specific conditions, the lasso's solutions can exactly recover the original $\ell_0$-norm solution \cite{candes2006robust,chen2001atomic}.

	The remaining part of this paper is organized as follows. In \Cref{sec.notation}, we introduce the notation and the definition of c-MIR cuts. In \Cref{sec.mw}, we revise the MW heuristic and its implementation, and we show its limitations in an example. In \Cref{sec.card}, we provide the cardinality optimization formulation of the aggregation,  we show that the MW aggregation heuristic is a stepwise selection heuristic, which is a relatively ancient approach in statistical community  \cite{draper1998applied,hastie2017extended}, and we introduce the lasso approximation of the cardinality optimization problem. In \Cref{sec.opt}, using the lasso approximation, we propose an aggregation algorithm that solves several linear programs (LPs) to pursue a sparse aggregated base inequality and project out integer variables. In \Cref{sec.exper}, we present the experiment results. In \Cref{sec.conclu}, we conclude this paper with the future research prospect. This perspective not only sheds light on the structure of the aggregation process but also provides a foundation for our further refinements.

	\section{Notation and c-MIR cuts}\label{sec.notation}
	In this section, we provide the necessary notation and the definition of c-MIR cuts \cite{marchand2001aggregation}.
	Let $[t_1:t_2]$ denote the sequence $(t_1,\dots, t_2)$, and let $[t]$ be an abbreviation for $[1:t]$.
	Denote $x^+ = \max(x, 0)$. Let $\norm{\cdot}_0$ and $\norm{\cdot}_1$ denote the $\ell_0$-norm and $\ell_1$-norm, respectively. Given a matrix $A \in \bR^{m \times n}$, $A_{i,j}$ is the entry at row index $i \in [m]$ and column index $j \in [n]$. To index its submatrices, $i$ and/or $j$ can be replaced with index sets $I \subseteq [m]$ and/or $J \subseteq [n]$. Additionally, when using a matrix index, the symbol $:$ denotes the full index set $[m]$ or $[n]$, depending on the context.

	Next, we define the c-MIR cut for the following mixed knapsack set:
	\begin{equation}
		\label{eq.mks}
		X = \left\{ (s,z) \in \bR_+ \times \bZ_+^q : \sum_{j \in [q]} a_j z_j \leq b + s, \; z_j \leq u_j \text{ for } j \in [q] \right\}.
	\end{equation}%\todo{KB: the intersection with $X'$ here doesn't make sense to me. Is it intended for this to be here? liding: it will used later}

	Let \( (T, U) \) be a partition of \( [q] \) and let \( \delta > 0 \). We can complement variables in \( U \) and scale the constraint $\sum_{j \in [q]} a_j z_j \leq b + s$  by dividing it by \( \delta \). This results in the single-row relaxation:
	\begin{equation}
		\label{eq.mks2}
		X' = \left\{ (s,z) \in \bR_+ \times \bZ_+^q : \sum_{j \in T} \frac{a_j}{\delta} z_j + \sum_{j \in U} \frac{-a_j}{\delta} (u_j - z_j) \leq  \frac{b - \sum_{j \in U} a_j u_j}{\delta} + \frac{s}{\delta} \right\}.
	\end{equation}
	Then, the c-MIR inequality derived from $X'$ is as follows:
	\begin{equation}
		\label{ineq.cmir}
		\sum_{j \in T} G(\frac{a_j}{\delta}) z_j + \sum_{j \in U} G(\frac{-a_j}{\delta}) (u_j - z_j) \leq \floor{\beta}  + \frac{s}{\delta (1-f)},
	\end{equation}
	where $
	\beta = (b - \sum_{j \in U} a_j u_j) / \delta$, $f = \beta - \lfloor \beta \rfloor$,
	and
	$
	G(d) = (\lfloor d \rfloor + (f_d - f)^+ / (1 - f)) \text{ with } f_d = d - \lfloor d \rfloor$.
	The c-MIR inequality is valid for $X'$ (hence for $X$), and it becomes the so-called MIR inequality when $U = \varnothing$  and $\delta = 1$ \cite{marchand2001aggregation,nemhauser1990recursive}.   The c-MIR inequalities can represent mixed cover inequalities, residual capacity inequalities, and weight inequalities.

	\section{The Marchand-Wolsey Heuristic} \label{sec.mw}
	In this section, we discuss the MW heuristic~\cite{marchand2001aggregation}, a procedure designed to generate c-MIR inequalities, and several enhancements of its implementations \cite{christophel2009separation,wolter2006implementation}

	% In the extreme case,   Benders decomposition projects out all continuous variables by enforcing $\lambda^\top A_y = 0$ ($A_y$ are columns associated with $y$).  In this sense, the resulting Benders cut $\lambda^\top Ax  \le \lambda^\top b $ is also an aggregated inequality.

	As we have mentioned, continuous variables are generally not preferred in the base inequality, since a base inequality consisting solely of continuous variables does not yield additional useful inequalities.  Although continuous variables can be included when deriving c-MIR cuts \eqref{ineq.cmir}, it is often beneficial to project out some continuous variables to enhance the effectiveness of the cuts. If we enforce the use of a row in the aggregated base inequality, since MILP instances are usually sparse, then there may  be only a few nonzero continuous variables in the enforced row to eliminate. Moreover, the complete elimination of continuous variables is not necessary. This characteristic marks a key distinction from Benders decomposition.

	The bound distance is a measure used in \cite{marchand2001aggregation,wolter2006implementation} to determine ``bad'' continuous variables that will be preferred for elimination. To illustrate this, we divide the linear constraints $Ax - b \le 0$  of the MILP \eqref{eq.milp} into two parts: the bound constraint part and the normal constraint part. For brevity, we denote $p+q$ by $n$.  The bound constraints, which include variable bounds and implied variable bounds, %\todo{KB: expanded the wording here to make it more specific}
	are expressed in the set
	\begin{equation}
		X_1  =  \{x  \in \bR^n ~:~ \forall j \in [n], x_j \le u_j, \text{ and } \forall j \in [p], j' \in J_j,   x_j \le u_{jj'} + d_{jj'} x_{j'}  \},
	\end{equation}
	where $J_j$ is a subset of $[n+1: n + q]$. The remaining normal constraints define the set $X_2 = \{x  \in \bR^n , A'x - b' \le 0\}$, where $A' = A_{\bar{I},:}$ and $b' = b_{\bar{I}}$ with $\bar{I} \subseteq I$ being the index set for the remaining constraints.  Thus, the MILP \eqref{eq.milp} reads as $\min_{x \in  \bR^p \times \bZ^{q}, x \in X_1 \cap X_2} c^\top x$. Given a point $\relx{x}$ to separate,  for $j \in [p]$, the bound distance of $\relx{x}_j$ is defined as
	\begin{equation}
		\label{eq.bddist}
		\bd(\relx{x}_j) =   \min(\{u_j\} \cup \{u_{jj'} + d_{jj'}\relx{x}_{j'} \}_{j' \in J_j}) - \relx{x}_j.
	\end{equation} In the heuristic, eliminating continuous variables with large bound distances is preferred, which we will explain later. A continuous variable with a nonzero bound distance is called a ``bad variable''.
	We next revise the three basic steps of the MW heuristic, which was presented at a high level in \cite{marchand2001aggregation}:
	\begin{enumerate}
		\item \textbf{Normal constraint aggregation}:  The MW heuristic identifies bad continuous variables, then it employs a stepwise selection algorithm to determine a few nonzero factors $\lambda_{\bar{I}}$. The aggregated inequality $\lambda_{\bar{I}}^\top(A'x - b') \le 0$ should be sparse, and some bad continuous variables are projected out.
		\item \textbf{Bound substitution}: %Some continuous variables
		Each remaining continuous variable $x_j$ in the aggregated inequality is substituted by an expression of its bound and a nonnegative continuous slack variable $y_j$: $x_j = u_j - y_j$ or $x_j = u_{jj'} + d_{jj'} z_{j'} - y_j$. %\todo{KB: I see what was intended in this paragraph, but strictly speaking, one can't substitute variables by inequalities. So I have changed this (and the description of $s$) to make this formally correct. I have also adjusted the explanation to fit the papers by MW and Wolter more closely.} % their bounds $y_j \le u_j$ or  $ y_j \le u_{jj'} + d_{jj'} z_{j'}$.
		%Thus, more continuous variables are projected out.
		This yields the inequality defining the mixed knapsack set $X$, where $s$ in \eqref{eq.mks} is taken as the weighted sum of all %unprojected
		continuous variables $y$, and $z$ in \eqref{eq.mks} contains all integer variables present in the base inequality.
		\item \textbf{Separation}: The parameters  \( (T,U) \) and \( \delta \) are determined, and the associated c-MIR inequality \eqref{ineq.cmir} is generated. 
	\end{enumerate}

\begin{figure}[ht]
		\centering
		\includegraphics[scale=0.4]{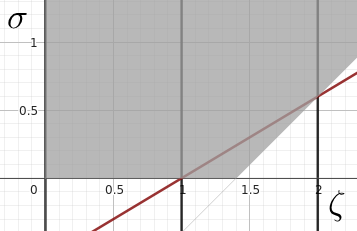}
		\caption{Example of a split cut}\label{fig:bound-distance-example}
	\end{figure}

	The MW heuristic's preference for eliminating continuous variables with large bound distances is explained as follows.
	We can reduce the set $X'$ in \eqref{eq.mks2}  to a simple two-variable mixed-integer set $\{(\sigma,\zeta) \in \bR_+ \times \bZ: \zeta \le \beta + \sigma \}$, and the c-MIR inequality reduces to the split inequality $\zeta  \le \floor{\beta} + \sigma/(1- (\beta - \floor{\beta})) $ valid for this set \cite{cornuejols2008valid,marchand2001aggregation}. See \Cref{fig:bound-distance-example} for an example. The gray area indicates the set of nonnegative points satisfying $\zeta \le \beta + \sigma$, and the vertical lines represent sets of integer feasible solutions. The line passing through (0,1) and (2,0.6) shows the boundary of the set of solutions satisfying the split inequality.

	Observe that feasible solutions of the LP relaxation that are cut off by the split inequality have relatively small values of $\sigma$. Further, the smaller the value of $\sigma$, the larger the corresponding interval that is cut off.  In the extreme case, consider a point $(\relx{\sigma}, \relx{\zeta})$. If $\relx{\sigma}$ is zero (\ie at its bound) and the base inequality $\zeta \le \beta + \sigma$ has zero slack (\ie $\relx{\zeta} =\relx{\beta}$), then the point is cut off by the split inequality. An aggregation that eliminates continuous variables with large bound distances, followed by the substitution of remaining continuous variables using their bounds, represents a shift into a space where the remaining nonnegative continuous variable ($\sigma$ in the example above) tends to be close to 0. This results in a higher likelihood that a strongly violated MIR inequality will be found.
	
	The MW heuristic differentiates between normal constraints and variable bound constraints. However, the second bound substitution step can be interpreted as a specialized
	constraint aggregation procedure applied to variable bound constraints, where the absolute values of the corresponding factors are restricted to either zero or one by the MW heuristic. We will show that, in \Cref{sec.opt}, our aggregation algorithm can unify the aggregation of two types of constraints, so the algorithm has no restrictions on factors for bounds.

	SCIP's implementation of the aggregation algorithm includes several enhancements, such as preprocessing and the use of row scores. Normal constraints' rows containing bad variables are identified as ``useful rows'' (that is, useful for eliminating bad variables). The preprocessing step reduces the LP's rows to these useful rows, effectively pruning the search space of $\lambda$.
	%This helps to prune the search space of $\lambda$.
	Row scores are used to select rows to combine. The score of a row is a weighted sum of its dual value,  density, slack, fractionalities of its integer variables, and ``fractionalities'' of its continuous variables defined through slacks of implied variable bounds. %\todo{KB: fractionalities of integer variables? Or are bound distances included in this as 'fractionalities' of continuous variables? liding: it is a bit complicated concept, and I added some illustration. KB: I elaborated a bit further.}

	\begin{algorithm}
		\caption{The preprocessing algorithm}
		\begin{algorithmic}[1]
			\Require the LP data $(A,b)$.
			\State compute bound distances of continuous variables;
			\State identify the index set $J$ of bad variables and the index set $I$ of useful rows;
			\State compute scores of useful rows;
			\State sort the indices of bad variables in decreasing order of bound distance;
			\State sort the indices of useful rows in decreasing order of score;
		\end{algorithmic}
		\label{algo.prep}
	\end{algorithm}

	The whole preprocessing algorithm is presented in \Cref{algo.prep}. Finally, we outline SCIP's implementation of the MW aggregation heuristic in \Cref{algo.mw}. %Note that the combination of aggregated rows with a new row in Line \ref{algo.try} may not be accepted in Line \ref{algo.accept}.

	%\todo{KB: does the starting row contain only bad variables? This doesn't look right. (I have replaced J with : now) liding: yes, only such rows, is J right here?}
	\begin{algorithm}
		\caption{The MW aggregation heuristic}
		\begin{algorithmic}[1]
			\Require relaxation solution $\relx{x}$, the LP data $(A,b)$, the starting row index $i_0$, and the maximum round MAXAGGR of aggregations.
			\State call the preprocessing \Cref{algo.prep};
			\State in the aggregated constraint $\alpha x - \beta \leq 0$, set $\alpha \leftarrow A_{i_0,:}$ and $\beta \leftarrow b_{i_0}$ %add a starting row $r_{i_0} = A_{i_0,:}x$ to the aggregated row $r=a^\top x$; %$r_{i_0} = A_{i_0,J}x_J  $ to the aggregated row $r=a^\top x_J$;
			\State set the aggregation counter $c \leftarrow  0$;
			\State set the set $I'$ of used rows  to $\{i_0\}$;
			\State set the set $J'$ of eliminated bad variables  empty;
			\For{each  sorted index $j \in J \smallsetminus J'$ of bad  variables} \label{restart}
			\For{each sorted index $i \in I \smallsetminus I'$ of useful rows for which $A_{ij} \ne 0$}
			\State find the factor $\lambda_{i}$ such that  the coefficient vector $\alpha^{\text{new}} = \alpha + \lambda_{i}  A_{i,:}$ %$\alpha = a+ \lambda_{i}  A_{i,J}$
			has $\alpha^{\text{new}}_j  = 0$; \label{algo.try}
			\If{there exists $j \in J'$ such that $\alpha^{\text{new}}_j \ne 0$} \label{badaggr} \label{algo.accept}
			\State \textbf{continue}~; \Comment{do not add rows that reintroduce eliminated bad variables}
			\EndIf
			\State update: $\alpha  \leftarrow  \alpha^{\text{new}}$, $\beta \leftarrow \beta + \lambda_ib_i$, $I' \leftarrow I' \cup \{i\}$, $J' \leftarrow J' \cup \{j\}$, and $c \leftarrow c + 1$;
			\State try bound substitution and separation on the aggregated constraint $\alpha^\top x  - \beta \le 0 $; %$a^\top x_J  - \lambda_{I'}^\top b_{I'}\le 0 $;
			\If{$c < $ MAXAGGR}
			\State  \textbf{break}~; \Comment{proceed to next bad variable}%\textbf{goto} Line~\ref{restart};
			\Else
			\State \Return;
			\EndIf
			\EndFor
			\EndFor
		\end{algorithmic}
		\label{algo.mw}
	\end{algorithm}

	The MW heuristic cannot adjust the factor $\lambda_i$ after the row $r_i$ is added into the aggregated row. Consequently, it cannot explore all possible combinations of rows. The following example illustrates a limitation of the MW heuristic, where it fails to eliminate all bad variables, even if there exists such an aggregation.

	\begin{example}
		\label{example.1}
		Consider the MILP's constraint $Ax \le b, x \ge 0$ as follows:
		\begin{subequations}
			\label{eq.example}
			\begin{align}
				1/3x_1 +  x_2 - 2/3x_3\le& 1, \label{eq.example1} \\
				2/3x_1 - 1/3 x_2 - 4/3 x_3 + x_4\le& 1, \label{eq.example2} \\
				-1/3 x_2 + x_3\le& 1, \label{eq.example3}
			\end{align}
		\end{subequations}
		where all rows are useful,  and $x_2$ and $x_3$ are bad variables. A vector $\lambda = (1,1,2)$ exists such that the corresponding aggregated inequality $x_1  + x_4\le 4$ eliminates all bad variables. However, the MW heuristic cannot find a projected inequality starting from any row. For example, starting from the first row \eqref{eq.example1},  assume that the MW heuristic chooses the second row \eqref{eq.example2} to eliminate $x_2$. In this case, the heuristic greedily select the factor $3$ for \eqref{eq.example2} and does not select any factor for \eqref{eq.example3}, then  we obtain the aggregated inequality $7/3x_1 -14/3 x_3 + 3x_4 \le 4$ using $\lambda=(1,3,0)$.
	\end{example}
	Hence, a systematic method is required to explore the space of $\lambda$, rather than relying on the stepwise activation and fixing of nonzero factors. In the following sections of this paper, we propose an optimization-based approach to address this need.
	\section{Cardinality Formulation of the Aggregation Problem} \label{sec.card}

	In this section, we formulate the aggregation problem as a cardinality optimization problem, more precisely a cosparsity problem \cite{nam2013cosparse}. We cast the problem as the standard weighted regularized cardinality optimization problem. We show that the MW aggregation heuristic is a stepwise selection heuristic, which was used to solve the $\ell_0$ minimization problem \cite{hastie2017extended}. Then, we present the lasso approximation of the cardinality optimization problem.

	We consider the following  cardinality optimization problem:

	\begin{equation}
		\label{eq.cosparsity}
		\min_{\lambda_I \in \bR^I_+, \lambda_{i_0} \ge 1}  \norm{ W \lambda_I^\top A_{I,J}}_0,
	\end{equation}
	where $A_{I,J}$ is the submatrix of $A$ with useful rows indexed by $I$ and bad variables indexed by $J$, $W = \Diag(w)$ is the diagonal weight matrix with diagonal entries  $w = (w_j)_{j \in J}= (\bd(\relx{x}_j))_{j \in J}$ (see \eqref{eq.bddist}). Note that we enforce the starting row $r_{i_0}$ by setting $\lambda_{i_0} \ge 1$ (otherwise, the minimum value is zero), and the objective is a weighted sum of nonzeros in the aggregated row: %\todo{KB: having the expanded expression is nice. Do we need the $\neq 0$ part, though? I see that this emphasizes that these are nonzeroes, but zeroes wouldn't contribute to the sum, anyway.}
	\begin{equation}
		\norm{ W \lambda_I^\top A_{I,J}}_0 =  \sum_{j \in J}w_j 1\{(\sum_{i \in I} A_{i,j} \lambda_i) \ne  0\}.
	\end{equation}
	The problem \eqref{eq.cosparsity} is a cosparsity problem, which presumes that the linear transformation $\lambda_I^\top A_{I,J} $ of $\lambda_I$ is sparse. The problem is nonconvex and is known to be $\mathcal{NP}$-hard \cite{tillmann2024cardinality}.

	Moreover, we want the aggregated base inequality to be sparse, so this translates to a requirement of sparse factors. Thus, we add $\ell_0$ regularization term on $\lambda_I$ to the objective, and we obtain the following regularized coparsity formulation of the aggregation problem:

	\begin{equation}
		\label{card.l0}
		\min_{\lambda_I \in \bR^I_+, \lambda_{i_0} \ge 1}  \norm{ W \lambda_I^\top A_{I,J}}_0 +  \norm{ W'  \lambda_I}_0,
	\end{equation}
	where  $W' = \Diag(w')$ is a diagonal weight matrix, and the diagonal entries  $w'_i$ ($i \in I$) is  a user-defined penalty.  For example, $w'_i$ could be the score of row $r_i$.

	We next cast the problem as the standard weighted regularized cardinality optimization problem. Let us introduce variable $\mu$ representing $\lambda_I^\top A_{I,J}$, and this leads to the weighted version of the $\ell_0$ minimization problem:
	\begin{equation}
		\label{l0}
		\min_{\lambda_I \in \bR^I_+, \lambda_{i_0} \ge 1, \lambda_I^\top A_{I,J} = \mu }  \norm{ W \mu}_0 +  \norm{ W'  \lambda_I}_0.
	\end{equation}

	Forward stepwise selection \cite{draper1998applied,hastie2017extended} and the lasso \cite{tibshirani1996regression,chen2001atomic} (aka $\ell_1$-norm approximation) are two popular methods for approximately solving the $\ell_0$ minimization problem \eqref{l0}. Forward stepwise selection can be seen as a ``greedy'' algorithm, updating the active set by one variable at each step, instead of re-optimizing over all possible subsets of a given size; in turn, the lasso can be seen as a more “democratic” version of forward stepwise, updating the coefficients to maintain equal absolute correlation of all active variables with the objective.

	The MW heuristic can be viewed as a forward stepwise selection algorithm for the $\ell_0$ formulation \eqref{card.l0}
	of the aggregation problem. In the statistical community, there are several known drawbacks of the forward stepwise selection heuristic. The number of the selected candidate variables may be smaller than the total number of final model variables \cite{hurvich1990impact}. Models that are created may be over-simplifications of the real models of the data \cite{roecker1991prediction}. These drawbacks imply that the MW heuristic cannot find complicated combinations of rows, which was already demonstrated in \Cref{example.1}.

	We adopt the lasso approach to tackle the cardinality formulation \eqref{card.l0}, %\todo{KB: I am not sure this is the right way to phrase this, since the lasso, by itself, is not a new method for solving these kinds of problems - only the use of it in this specific context is new. liding: ok, remove new}
	and this yields our new heuristic to solve the aggregation problem. As  the  $\ell_1$-norm promotes sparsity, the lasso  approximation uses the $\ell_1$-norm to approximate the $\ell_0$-norm in \eqref{l0}:
	\begin{equation}
		\label{l1}
		\min_{\lambda_I \in \bR^I_+, \lambda_{i_0} \ge 1, \lambda_I^\top A_{I,J} = \mu }  \norm{ W \mu}_1 +  \norm{ W'  \lambda_I}_1.
	\end{equation}

	In our setting,  we let the penalty $w'_i$ be the nonnegative slack $ b_i - A_{i,:}^\top \relx{x}$ of the row. This penalty and $\ell_1$ approximation promote zero slacks, \ie tight rows.   This casts $\norm{ W'  \lambda_I}_1$ as $\sum_{i \in I}\lambda_i (b_i - A_{i,:}^\top \relx{x})$, and the  lasso approximation \eqref{l1} becomes:
	\begin{equation}
		\label{eq.lasso}
		\min_{\lambda_I \in \bR^I_+, \lambda_{i_0} \ge 1, \lambda_I^\top A_{I,J} = \mu }  \norm{ W \mu}_1 + \sum_{i \in I}\lambda_i (b_i - A_{i,:}^\top \relx{x}).
	\end{equation}
	Note that the term $\sum_{i \in I}\lambda_i (b_i - A_{i,:}^\top \relx{x})$ is precisely the slack of the aggregated base inequality.%\todo{KB: what about the non-bad variables? Do they not contribute to the slack? liding: that was a mistake, as the new does contain non-bad variables.}

	Thus, viewing instead $\norm{ W \lambda_I^\top A_{I,J}}_1$ as the regularization term, we reach another interpretation of the lasso approximation: the $\ell_1$ approximation aims at finding a sparse aggregated base inequality with a small slack (see \Cref{sec.mw} for the intuition). The lasso is a convex optimization problem (LP), so we can solve it to optimality.  Note that we can convert an equality row $\alpha^\top x = \beta $ into two inequality rows $\alpha^\top x - \beta \le 0 $ and $-\alpha^\top x + \beta \le 0  $ both of zero slacks. An optimal basic solution to \eqref{algo.lp} usually has at most one non-zero factor for these two rows.

	\section{An Optimization-Driven Aggregation Algorithm} \label{sec.opt}

	In this section, we present an LP-based algorithm to solve the aggregation problem. The algorithm solves several LPs to promote the sparsity of the solution.

	Our aggregation algorithm can unify the normal constraint aggregation and bound substitution steps of the MW heuristic. To this end, we modify the preprocessing algorithm by treating variable bound constraints as normal constraints, \ie useful rows can include variable bounds. Thus, the index set $I$ of useful rows is a subset of $[m]$ (not $\bar{I}$). This exploits more possibilities to project out bad variables. We next present the main aggregation algorithm.

	In the first factor search stage, our algorithm solves the LP \eqref{eq.lasso} to find a $\lambda_I$.  Denote by  $I' = \{i \in I: \lambda_i > 0\}$ the index set of active factors. Then, the aggregated base inequality is $\lambda_{I'} A_{I',J} x  - \lambda_{I'} b_{I'} \le 0 $. As lasso is sparse inducing, the nonzero factors in the solution  $\lambda_I$ are usually quite sparse. This implies that the  aggregated base inequality is sparse.   Moreover, the sparsity of $\lambda_{I'}^\top A_{I',J}$ is also promoted by the LP, and some bad variables are eliminated in the aggregated base inequality.

	If  $ \lambda_{I'}^\top A_{I',J}$ is not sparse enough, \ie the density of bad variables exceeds a threshold, then our algorithm moves to the second stage, wherein it applies iterative reweighting. The iterative reweighted algorithm \cite{candes2008enhancing} performs exceptionally well in locating sparse solutions to $\ell_0$ minimization, and it even outperforms $\ell_1$
	approximation. The key feature of the reweighted algorithm
	is the solution of a series of reweighted $\ell_1$
	approximation problems. In this stage,
	we use the iterative reweighted algorithm to construct and solve several modified LPs to find a sparser $\lambda$ supported on $I'$. As the active set $I'$ is already known, we are only interested in projecting out the bad variables, so our new LP has only the second term in its objective:

	\begin{equation}
		\label{lp.iter}
		\min_{\lambda_{I'} \in \bR^I_+, \lambda_{i_0} \ge 1}  \norm{ W \lambda_{I'}^\top A_{I',J}}_1.
	\end{equation}

	In each round of iterative reweighting, the algorithm dynamically adjusts the weight $W$ and resolves the modified LP \eqref{lp.iter}. The weights are updated with the use of the previous solution $\lambda$ to \eqref{lp.iter}:

	\begin{equation}
		\label{eq.reweight}
		w_i = \begin{cases}
			w_i / (\epsilon + \lambda_i),  &\lambda_i \ne 0,\\
			0,  & \lambda_i = 0,
		\end{cases}
	\end{equation}
	where $\epsilon$ is a small positive constant introduced to avoid numeric overflow. We solve the iteratively reweighted LPs \eqref{lp.iter} until $ \lambda_{I'}^\top A_{I',J}$ is sparse enough.

	We next present our sparsity-driven  LP-based aggregation algorithm in \Cref{algo.lp}.

	\begin{algorithm}
		\caption{Sparsity-driven  LP-based aggregation algorithm}
		\begin{algorithmic}[1]
			\Require relaxation solution $\relx{x}$, the LP data $(A,b)$, the starting row index $i_0$,  the maximum round MAXAGGR of aggregations, and the density threshold $\eta$.
			\State call the preprocessing \Cref{algo.prep};
			\State solve the lasso approximation \eqref{eq.lasso};
			\State set the set $I'$ of used rows  as support of nonzeros in $\lambda_I$;
			\State set the aggregation counter $c \leftarrow  0$;
			\While{$c$ is less than MAXAGGR}
			\State let $a = \lambda_{I'}^\top A_{I',:
			}$;
			\State try separation on the aggregated base inequality $a^\top x  - \lambda_{I'}^\top b_{I'
			}\le 0 $;
			\If{the density $|\{j \in J:a_j \ne 0\}| / |J|$ of bad variables is greater than the threshold $\eta$}
			\State  let $c \leftarrow c + 1$;
			\State update the weight $w$ as in \eqref{eq.reweight};
			\State solve the reweighted LP \eqref{lp.iter};
			\Else
			\State \Return;
			\EndIf
			\EndWhile
		\end{algorithmic}
		\label{algo.lp}
	\end{algorithm}

	\begin{example}
		We revise \Cref{example.1} and apply the LP-based aggregation algorithm instead of the MW heuristic. Assume that the value of the point $(x_1,x_2,x_3)$ in \eqref{eq.example} is zero, and the slacks of all constraints in \eqref{eq.example}  are zero. Recall that $x_2$ and $x_3$ are bad variables. We can express the corresponding lasso approximation as:
		\begin{equation}
			\min_{\lambda \in \bR^3_+, \lambda_{i_0} \ge 1} w_2|\lambda_1 - 1/3 \lambda_2  - 1/3 \lambda_3 | + w_3 |-2/3\lambda_1 - 4/3 \lambda_2  + \lambda_3 |,
		\end{equation}
		where $w_2$ and $w_3$ are non-negative weights (set to be equal to bound distances), and $i_0 \in \{1,2,3\}$. Note that the term $\sum_{i \in I}\lambda^\top_i (b_i - A_{i,J}^\top \relx{x}_J)$  is null due to the slacks being zero. The solution $(1,1,2)$ is optimal with objective value zero regardless of values of $w_2,w_3$, and this solution eliminates all bad variables. Thus, the LP approach finds a suitable aggregation, and it is stable under the perturbation of weights.
	\end{example}

	\section{Experimental Comparison} \label{sec.exper}
	In this section, we evaluate the performance of SCIP on MIP problems and compare it with the performance achieved with the addition of our new aggregation-based cut separator.
	\subsection{Implementation Details}
	The \Cref{algo.mw}  and \Cref{algo.lp} have been implemented in the open-source MIP solver SCIP \cite{bolusani2024scipoptimizationsuite90}.  At each node of the branch-and-bound tree, several rounds of LPs are solved, and separators are called to strengthen the LP relaxations. SCIP has the default aggregation cut separator based on the MW heuristic, and we implement the new aggregation cut separator based on  \Cref{algo.lp}.
	In a single separation round, the cut separators generate multiple aggregations using different starting rows and try separations on the aggregated base inequalities. The generation of flow cover inequalities \cite{gu1999lifted}, knapsack cover inequalities \cite{letchford2019lifted}, and c-MIR inequalities on an aggregated base inequality is tried, with the first two types of inequalities being concrete models of c-MIR inequalities.

	There are two main differences between our aggregation-based cut separator and the default one, which result in our separator generating fewer cuts:
	\begin{itemize}
		\item The first difference is by design to prevent the generation of repetitive cuts. Different starting rows in \Cref{algo.mw,algo.lp} may lead to the same aggregated base inequality. To avoid redundancy and conserve computational resources, our separator only utilizes a starting row if it has not been used in previous aggregations.
		\item The second difference does not directly involve aggregation. Unlike our method, the default separator attempts to separate cuts on every row of a MILP instance, even if the instance is purely an ILP. For example, the default separator might separate cuts for a knapsack constraint.
	\end{itemize}
	The two separators share the same parameter settings, defined as follows along with their respective values:
	\begin{itemize}
		\item SEPA\_FREQ (10): the frequency defines the depth levels at which the separation is called. A frequency of 10 means that the separation is executed in depth 0, 10, 20, \dots of the branching tree.
		\item MAXAGGR (6 at root node, 3 at non-root nodes): the maximum round of aggregations, see \Cref{algo.mw,algo.lp}.
	\end{itemize}

	For the cut separator based on our new \Cref{algo.lp}, we limit the numbers of bad variables and useful rows to 50 and 5000, respectively. The aggregation LPs are solved using the primal simplex algorithm of CPLEX, so that our algorithm can benefit from warm-start after changes of objectives and lower bounds of factor variables during a separation round.

	\subsection{Results}
	All experiments were conducted on a cluster equipped with Intel Xeon Gold 5122 CPUs, operating at 3.60 GHz with 96 GB of RAM. Each test run was executed on a single thread with a time limit of two hours. To ensure accurate comparisons despite performance variability~\cite{variability}, we used a large and diverse test set, specifically the MIPLIB 2017 benchmark set~\cite{miplib2017}. We used random seeds ranging from 0 to 4, resulting in a total of 1200 instance-seed pairs.

	Table~\ref{tbl:MIP5seeds} compares the performance of default SCIP (``SCIP-default''), with SCIP using the new aggregation-based cut separator (``SCIP+aggr.''). It shows the numbers of solved instances, the shifted geometric means of running time, and the number of branch-and-bound tree nodes. In the computation of geometric means, we use a shift of 1 second and 100 nodes, respectively. Relative differences in time and nodes for each category are shown in the last two columns of the table, with relative differences that exceed 1\% highlighted.

	We present these values for various subsets of the test instances. The category ``all'' encompasses all 1200 test instances. ``solved-by-both'' refers to instances that were solved to optimality by both configurations. The ``affected'' category highlights instances where performance varied between the two configurations, indicated by a difference in the total number of LP iterations. Further, we compare the performance on different categories of affected instances: instances that both configurations solved in under 100 seconds (``solved-under-100s''), instances that took more than 100 seconds to solve by either configuration (``solved-over-100s''), and instances that only one of the configurations managed to solve (``solved-by-one'').

	\begin{table}[h]
		\caption{Performance Comparison of SCIP-default with SCIP+aggr.}
		\label{tbl:MIP5seeds}
		\footnotesize
		\begin{tabularx}{\textwidth}{@{}l@{\;\;\extracolsep{\fill}}rrrrrrrrr@{}}
			\toprule
			&           & \multicolumn{3}{c}{SCIP-default} & \multicolumn{3}{c}{SCIP+aggr} & \multicolumn{2}{c}{relative} \\
			\cmidrule{3-5} \cmidrule{6-8} \cmidrule{9-10}
			Subset              & instances & solved & time & nodes & solved & time & nodes & time & nodes \\
			\midrule
			all & 1200 & 714 & 936.9 & 6436 & \win{716} & 939.9 & 6151 & 1.00 &  - \\
			solved-by-both & 706 & 706 & 226.4 & 2711 & 706 & 228.6 & 2584 & 1.01 & \win{0.95} \\
			\cmidrule{1-10}
			affected & 330 & 320 & 287.9 & 5135 & 322 & 289.9 & 4629 & 1.01 & - \\
			solved-under-100s & 98 & 98 & 17.7 & 375 & 98 & 20.9 & 320 & \loss{1.18}& \win{0.85}  \\
			solved-over-100s & 214 & 214 & 789.4 & 13133 & 214 & 751.6 & 11931 & \win{0.95} & \win{0.91}  \\
			solved-by-one & 18 & 8 & 5476.5 & 39786 & 10 & 4719.1 & 37720 & \win{0.86}  & - \\
			\bottomrule
		\end{tabularx}
	\end{table}

	Across the entire data set, there is minimal difference in time or number of solved instances, with only two additional instances solved by SCIP+aggr. In the category ``solved-by-both'', we observe a 5\% decrease in the number of nodes.
	In the category ``solved-under-100s'', there is a notable 17\% increase in time, despite a 15\% reduction in the number of nodes. This increase is expected, as the separation algorithm is expensive for easy-to-solve problems. In contrast, in the category ``solved-over-100s,'' we observed a 5\% decrease in time and a 9\% reduction in nodes, which indicates that our cut separator helps solve harder problems more efficiently.
	To assess the statistical significance in the reduction in the number of nodes from ``SCIP-default'' to ``SCIP+aggr'', we performed a version of the Wilcoxon signed-rank test~\cite{berthold2015heuristic}. The test, applied to instances solved by both models, yielded a p-value below 0.01, indicating a significant difference in the number of nodes between the two settings.
	% \todo{KB: what about the subsets? Or you didn't want to go into that much detail? I guess it's fine like this (so feel free to comment this out).}
	For the limited subset of affected instances ``solved-by-one'' SCIP+aggr was 14\% faster.

	Finally, to complement our performance analysis, we examine the sparsity of the aggregated rows produced by \Cref{algo.mw} and \Cref{algo.lp}. In \Cref{tab:metrics}, we report the average number of bad columns in an aggregated row, the average of the total number of distinct bad columns across all rows used for aggregation, their ratio, and the average number of rows used for aggregation.

	The results indicate that rows aggregated by \Cref{algo.lp} contain significantly fewer bad columns on average—only 0.37 compared to 2.45 by \Cref{algo.mw}.
	%Furthermore, our approach projects out a higher percentage of bad columns; specifically, only
	This corresponds to 20\% of bad columns remaining after aggregation compared to 52\% with the MW heuristic.
	Additionally, we consider the number of rows used in the aggregation, which typically implies the sparsity of the resulting aggregated row, as its number of nonzeros is upper-bounded by the total number of nonzeros of the used rows. The statistics show that our method utilizes fewer rows in aggregation. This reveals that our aggregation process naturally produces sparse cuts, which is a desirable property.
	Both aggregation algorithms utilize only a small subset of rows for aggregation, with the number  being significantly smaller than the total number of useful rows (note that our separator imposes a limit of 5000 useful rows). This behavior may be linked to the inherent sparsity of instances in the MIPLIB 2017 benchmark. This suggests that one may obtain a potential performance improvement using LP solvers that can effectively exploit sparsity.

	All the numbers in \Cref{tab:metrics} are small, but the size of lasso approximation \eqref{eq.lasso} can be much larger than these numbers. This indicates that important (non-zero) factors are confined to a low-dimensional subspace of $\bR^I$. Similar scenarios arise in fields such as signal processing, machine learning, and statistics. Notably, the lasso method has been demonstrated \cite{candes2005decoding} to effectively recover such low-dimensional subspaces, given the sparsity of the model. This further justifies our approach.

	\begin{table}
		\centering
		\caption{Sparsity of rows aggregated by \Cref{algo.mw,algo.lp}}
		\label{tab:metrics}
		\begin{tabular}{lrr}
			\toprule
			&       \Cref{algo.mw} & \Cref{algo.lp} \\
			\midrule
			bad-cols &  2.45 &  0.37 \\
			total-bad-cols &    4.72 &  1.87 \\
			ratio &     0.52 &  0.20 \\
			used-rows &     3.10 &  2.24 \\
			\bottomrule
		\end{tabular}
	\end{table}

	\section{Conclusion} \label{sec.conclu}

	We developed a new approach for aggregating MILP constraints to generate mixed-integer knapsack sets for which strong c-MIR cutting planes can be constructed.
	Similarly to the standard greedy heuristic, our algorithm aims to project out as many bad continuous variables as possible within given working limits, the ``badness'' being defined as slacks of variable bound or implied variable bound constraints.
	The algorithm approximately solves the optimization problem of finding such aggregation weights that minimize the number of non-zeroes corresponding to bad continuous variables.
	This problem is formulated as an $\ell_0$ norm minimization problem and solved with the use of the lasso technique.
	If the number of bad variables in the aggregated row still exceeds a threshold, iterative reweighting is performed to project out further bad variables.

	We give examples to 1) illustrate the connection between the continuous variable bound distance and the potential for generating a strong rounding cut for a given mixed-integer set, and 2) show that our approach has a greater capacity to project out bad variables than the standard heuristic.
	A unifying perspective on the standard and new aggregation algorithms is presented, demonstrating that both represent different approaches to approximately solving the $\ell_0$ norm minimization problem: in particular, the standard algorithm is a stepwise selection heuristic for this problem.

	The algorithm is implemented within the constraint integer programming solver SCIP and tested by running SCIP with the standard and the new algorithms on the MIPLIB 2017 test set.
	When considering averages over all instances, we observe only a very small impact on mean performance, but dividing the instances into subsets of easier instances where SCIP solves the problem in under 100s with both configurations compared, and challenging instances which take at least 100s with at least one of the configurations, reveals a connection between instance difficulty and performance impact.
	While the new algorithm results in an 18\% increase in time on easier instances, likely due to the increased computational cost of aggregation, it yields a 5\% decrease in mean solver running time on challenging instances.
	This indicates that a flexible approach to select the right aggregation algorithm based on estimated instance difficulty could leverage the strength of both methods, and could be a subject of future work.

	Further work will, on the one hand, delve deeper into the theory of the connection between the properties of mixed-knapsack sets and the quality of c-MIR cutting planes constructed for these sets, and perform a theoretical analysis of the behavior of aggregation algorithms and the quality of the formulations they produce.
	On the other hand, it will extend these ideas to the mixed-integer nonlinear programming case, starting with conic MIPs which, on the one hand, have a well-defined structure, and, on the other hand, encompass a rich variety of problems beyond MILPs.  Further, minimization of the sum of absolute values in \eqref{eq.lasso} is a classical problem in nonsmooth convex  optimization. Due to this, we want to find a more efficient algorithm than the primal simplex for faster aggregation, and a possible candidate is the primal-dual hybrid gradient descent algorithm \cite{chambolle2011first}. A more open problem is how to avoid generating repetitive aggregations.
	
	\ifthenelse{\boolean{esa}}
	{
		
	}{
		\section*{Acknowledgements}
		We want to thank Felipe Serrano for his insightful comments on this paper. The work for this article has been conducted in the Research Campus MODAL funded by the German Federal Ministry of Education and Research (BMBF) (fund numbers 05M14ZAM, 05M20ZBM).
	}

	%MINLP (conic MIPs). theory for the heuristic, on the bound distance. Maybe mention column generation and gradient descent?
	
	\bibliography{reference}

\begin{thebibliography}{10}

\bibitem{achterberg2007constraint}
Tobias Achterberg.
\newblock Constraint integer programming.
\newblock 2007.

\bibitem{achterberg2010mcf}
Tobias Achterberg and Christian Raack.
\newblock The {MCF}-separator: detecting and exploiting multi-commodity flow
  structures in {MIP}s.
\newblock {\em Mathematical Programming Computation}, 2(2):125--165, 2010.

\bibitem{atamturk2002capacitated}
Alper Atamt{\"u}rk.
\newblock On capacitated network design cut--set polyhedra.
\newblock {\em Mathematical Programming}, 92:425--437, 2002.

\bibitem{balas1979disjunctive}
Egon Balas.
\newblock Disjunctive programming.
\newblock {\em Annals of discrete mathematics}, 5:3--51, 1979.

\bibitem{bnnobrs1962partitioning}
Jacques~F Benders.
\newblock Partitioning procedures for solving mixed-variables programming
  problems.
\newblock {\em Numer. Math}, 4(1):238--252, 1962.

\bibitem{berthold2015heuristic}
Timo Berthold.
\newblock {\em Heuristic algorithms in global {MINLP} solvers}.
\newblock Verlag Dr. Hut, 2015.

\bibitem{bienstock1996capacitated}
Daniel Bienstock and Oktay G{\"u}nl{\"u}k.
\newblock Capacitated network design—polyhedral structure and computation.
\newblock {\em Informs journal on Computing}, 8(3):243--259, 1996.

\bibitem{bolusani2024scipoptimizationsuite90}
Suresh Bolusani, Mathieu Besançon, Ksenia Bestuzheva, Antonia Chmiela, João
  Dionísio, Tim Donkiewicz, Jasper van Doornmalen, Leon Eifler, Mohammed
  Ghannam, Ambros Gleixner, Christoph Graczyk, Katrin Halbig, Ivo Hedtke,
  Alexander Hoen, Christopher Hojny, Rolf van~der Hulst, Dominik Kamp, Thorsten
  Koch, Kevin Kofler, Jurgen Lentz, Julian Manns, Gioni Mexi, Erik Mühmer,
  Marc~E. Pfetsch, Franziska Schlösser, Felipe Serrano, Yuji Shinano, Mark
  Turner, Stefan Vigerske, Dieter Weninger, and Liding Xu.
\newblock The {SCIP} {O}ptimization {S}uite 9.0, 2024.
\newblock URL: \url{https://arxiv.org/abs/2402.17702}, \href
  {https://arxiv.org/abs/2402.17702} {\path{arXiv:2402.17702}}.

\bibitem{bonami2008projected}
Pierre Bonami, G{\'e}rard Cornu{\'e}jols, Sanjeeb Dash, Matteo Fischetti, and
  Andrea Lodi.
\newblock Projected {C}hv{\'a}tal--{G}omory cuts for mixed integer linear
  programs.
\newblock {\em Mathematical Programming}, 113(2):241--257, 2008.

\bibitem{candes2006robust}
Emmanuel~J Cand{\`e}s, Justin Romberg, and Terence Tao.
\newblock Robust uncertainty principles: Exact signal reconstruction from
  highly incomplete frequency information.
\newblock {\em IEEE Transactions on information theory}, 52(2):489--509, 2006.

\bibitem{candes2005decoding}
Emmanuel~J Candes and Terence Tao.
\newblock Decoding by linear programming.
\newblock {\em IEEE transactions on information theory}, 51(12):4203--4215,
  2005.

\bibitem{candes2008enhancing}
Emmanuel~J Candes, Michael~B Wakin, and Stephen~P Boyd.
\newblock Enhancing sparsity by reweighted $\ell_1$ minimization.
\newblock {\em Journal of Fourier analysis and applications}, 14:877--905,
  2008.

\bibitem{caprara19960}
Alberto Caprara and Matteo Fischetti.
\newblock $\{$0, 1/2$\}$-{C}hv{\'a}tal-{G}omory cuts.
\newblock {\em Mathematical Programming}, 74:221--235, 1996.

\bibitem{caprara2000separation}
Alberto Caprara, Matteo Fischetti, and Adam~N Letchford.
\newblock On the separation of maximally violated mod-k cuts.
\newblock {\em Mathematical Programming}, 87:37--56, 2000.

\bibitem{chambolle2011first}
Antonin Chambolle and Thomas Pock.
\newblock A first-order primal-dual algorithm for convex problems with
  applications to imaging.
\newblock {\em Journal of mathematical imaging and vision}, 40:120--145, 2011.

\bibitem{chen2001atomic}
Scott~Shaobing Chen, David~L Donoho, and Michael~A Saunders.
\newblock Atomic decomposition by basis pursuit.
\newblock {\em SIAM review}, 43(1):129--159, 2001.

\bibitem{christophel2009separation}
Philipp~M Christophel.
\newblock {\em Separation algorithms for cutting planes based on mixed integer
  row relaxations}.
\newblock PhD thesis, Paderborn, Univ., Diss., 2009, 2009.

\bibitem{chvatal1973edmonds}
Vasek Chv{\'a}tal.
\newblock Edmonds polytopes and a hierarchy of combinatorial problems.
\newblock {\em Discrete mathematics}, 4(4):305--337, 1973.

\bibitem{cook1990chvatal}
William Cook, Ravindran Kannan, and Alexander Schrijver.
\newblock Chv{\'a}tal closures for mixed integer programming problems.
\newblock {\em Mathematical Programming}, 47(1):155--174, 1990.

\bibitem{cornuejols2008valid}
G{\'e}rard Cornu{\'e}jols.
\newblock Valid inequalities for mixed integer linear programs.
\newblock {\em Mathematical programming}, 112(1):3--44, 2008.

\bibitem{draper1998applied}
Norman~R Draper and Harry Smith.
\newblock {\em Applied regression analysis}.
\newblock McGraw-Hill. Inc, 1998.

\bibitem{fischetti2007optimizing}
Matteo Fischetti and Andrea Lodi.
\newblock Optimizing over the first {C}hv{\'a}tal closure.
\newblock {\em Mathematical Programming}, 110(1):3--20, 2007.

\bibitem{miplib2017}
Ambros Gleixner, Gregor Hendel, Gerald Gamrath, Tobias Achterberg, Michael
  Bastubbe, Timo Berthold, Philipp Christophel, Kati Jarck, Thorsten Koch, Jeff
  Linderoth, et~al.
\newblock {MIPLIB} 2017: data-driven compilation of the 6th mixed-integer
  programming library.
\newblock {\em Mathematical Programming Computation}, 13(3):443--490, 2021.

\bibitem{gomory1960algorithm}
Ralph~Edward Gomory.
\newblock An algorithm for the mixed integer problem.
\newblock {\em Report No. P-1885, The Rand Corporation, Santa Monica, CA.},
  1960.

\bibitem{gonccalves2005implementation}
Joao~PM Gon{\c{c}}alves and Laszlo Ladanyi.
\newblock An implementation of a separation procedure for mixed integer
  rounding inequalities.
\newblock {\em IBM Res. Report RC23686 (W0508-022), IBM}, 2005.

\bibitem{gu1999lifted}
Zonghao Gu, George~L Nemhauser, and Martin~WP Savelsbergh.
\newblock Lifted flow cover inequalities for mixed 0-1 integer programs.
\newblock {\em Mathematical Programming}, 85:439--467, 1999.

\bibitem{hastie2017extended}
Trevor Hastie, Robert Tibshirani, and Ryan~J Tibshirani.
\newblock Extended comparisons of best subset selection, forward stepwise
  selection, and the lasso.
\newblock {\em arXiv preprint arXiv:1707.08692}, 2017.

\bibitem{hurvich1990impact}
Clifford~M Hurvich and Chih—Ling Tsai.
\newblock The impact of model selection on inference in linear regression.
\newblock {\em The American Statistician}, 44(3):214--217, 1990.

\bibitem{letchford2019lifted}
Adam~N Letchford and Georgia Souli.
\newblock On lifted cover inequalities: A new lifting procedure with unusual
  properties.
\newblock {\em Operations Research Letters}, 47(2):83--87, 2019.

\bibitem{variability}
Andrea Lodi and Andrea Tramontani.
\newblock Performance variability in mixed-integer programming.
\newblock In {\em Theory driven by influential applications}, pages 1--12.
  INFORMS, 2013.

\bibitem{marchand19990}
Hugues Marchand and Laurence~A Wolsey.
\newblock The 0-1 knapsack problem with a single continuous variable.
\newblock {\em Mathematical Programming}, 85:15--33, 1999.

\bibitem{marchand2001aggregation}
Hugues Marchand and Laurence~A Wolsey.
\newblock Aggregation and mixed integer rounding to solve {MIP}s.
\newblock {\em Operations research}, 49(3):363--371, 2001.

\bibitem{mexi2024cut}
Gioni Mexi, Felipe Serrano, Timo Berthold, Ambros Gleixner, and Jakob
  Nordstr{\"o}m.
\newblock Cut-based conflict analysis in mixed integer programming.
\newblock {\em arXiv preprint arXiv:2410.15110}, 2024.

\bibitem{nam2013cosparse}
Sangnam Nam, Mike~E Davies, Michael Elad, and R{\'e}mi Gribonval.
\newblock The cosparse analysis model and algorithms.
\newblock {\em Applied and Computational Harmonic Analysis}, 34(1):30--56,
  2013.

\bibitem{nemhauser1990recursive}
George~L Nemhauser and Laurence~A Wolsey.
\newblock A recursive procedure to generate all cuts for 0--1 mixed integer
  programs.
\newblock {\em Mathematical Programming}, 46(1):379--390, 1990.

\bibitem{roecker1991prediction}
Ellen~B Roecker.
\newblock Prediction error and its estimation for subset-selected models.
\newblock {\em Technometrics}, 33(4):459--468, 1991.

\bibitem{tibshirani1996regression}
Robert Tibshirani.
\newblock Regression shrinkage and selection via the lasso.
\newblock {\em Journal of the Royal Statistical Society Series B: Statistical
  Methodology}, 58(1):267--288, 1996.

\bibitem{tillmann2024cardinality}
Andreas~M Tillmann, Daniel Bienstock, Andrea Lodi, and Alexandra Schwartz.
\newblock Cardinality minimization, constraints, and regularization: a survey.
\newblock {\em SIAM Review}, 66(3):403--477, 2024.

\bibitem{wolter2006implementation}
Kati Wolter.
\newblock Implementation of cutting plane separators for mixed integer
  programs.
\newblock {\em Dipolma thesis, Technische Universit{\"a}t Berlin}, 2006.

\end{thebibliography}
	
\end{document}